\documentclass[12pt,bezier]{article}
\usepackage{amsmath,amssymb,amsfonts,epsf,euscript}
\usepackage[usenames]{color}
\usepackage{graphicx}
\usepackage{epstopdf}
\usepackage{subfigure}
\usepackage{color}
\newcommand{\be}{\begin{equation}}
\newcommand{\ee}{\end{equation}}
\newtheorem{theorem}{Theorem}[section]
\newtheorem{lemma}[theorem]{Lemma}
\newtheorem{corollary}{Corollary}[section]
\newtheorem{example}{Example}[section]
\newtheorem{definition}{Definition}[section]
\newtheorem{remark}{Remark}[section]
\newtheorem{proposition}[theorem]{Proposition}
\renewcommand{\theequation}{\arabic{section}.\arabic{subsection}.\arabic{equation}}
\title{\bf\Large Structure of continuous-time ARMA process driven by semi-Levy measure}
\vspace{-1cm}
\author
{N. Modarresi\thanks{\scriptsize
Department of Mathematics and computer science, Allameh Tabataba'i University, Tehran, Iran.
E-mail:  n.modarresi@atu.ac.ir(N. Modarresi).}
\and  S. Rezakhah$^+$ \and S. Shoaee$^{+}$ \thanks{
Faculty of Mathematics and Computer Science, Amirkabir University of
Technology, 424 Hafez Avenue, Tehran 15914, Iran. E-mail:  rezakhah@aut.ac.ir(S. Rezakhah).
\hspace{3cm}
$^{++}$Department of Statistics, Faculty of Mathematical Sciences, Shahid Beheshti University, Tehran, Iran.
E-mail:  shirin\_shoaee@aut.ac.ir(S. Shoaee)
}}

\date{}
\begin{document}
\maketitle

\begin{abstract}
A class of continuous-time autoregressive moving average (CARMA) process driven by simple semi-Levy measure is defined and its properties are studied.
We discuss some new insights on the structure of the semi-Levy measure which is described as periodically divisible measure.
This consideration enable us to provide statistical property of the introduced process.
We show that this process is well defined without having to assume further conditions on the measure.
We find a kernel representation of the process and present the properties of first and second moments of it.
Finally we show the efficiency of our model by implying simulated data.\\

{\it AMS 2010 Subject Classification:} 60E07, 60G18, 60G51.\\

{\it Keywords:} Continuous time ARMA; Periodic random measure; Semi-Levy process.
\end{abstract}

\section{Introduction}
Continuous-time models for time series exhibit both heavy-tailed and long-memory behavior. Such models are of considerable interest, specially for the modeling of financial time series.
Early papers have studied the statistical analysis of continuous-time autoregressive (CAR) processes and continuous-time autoregressive moving average (CARMA) processes \cite{bro2}, \cite{bro3}, \cite{bro4}. Continuous-time models have also been utilized and analyzed successfully for the modeling of irregularly spaced data.

For the first time, Brockwell \cite{bro1} introduced the linear continuous-time model which is particularly advantageous for dealing with irregularly spaced data as continuous-time threshold ARMA$(p, q)$ process with $0 \leqslant q < p$. It provides the weak solution of a certain stochastic differential equation which is unique. Stramer et al. \cite{s3} investigated the existence and stability properties of these processes.
Properties of linear CARMA processes driven by second order Levy processes are examined and extended to include heavier tailed series which frequently encountered in financial applications \cite{bro1-1}.
Discrete time representations for data generated by a CARMA system with mixed stock and flow data are derived by Chambers et al. \cite{c0}.

Using the kernel representation of a Levy-driven CARMA process, the class of non-negative Levy-driven are extended by many authors for the non-monotone auto covariance functions. A class of fractionally integrated processes and also asymptotic properties of the CARMA processes are studied \cite{bro1-2}.
The second order Levy-driven CARMA models and some of their financial applications in particular to the modeling of stochastic volatility are discussed by Brockwell \cite{bro2}. Compound Poisson process and Levy-driven stationary Ornstein Uhlenbeck (OU) process and some examples of theses models are studied in \cite {b0}.
CARMA processes with a nonnegative kernel driven by a nondecreasing Levy process constituted a very general class of stationary nonnegative continuous time processes. The advantage of the nonnegativity of the increments of the driving Levy process is taken to develop a highly efficient estimation procedure for the parameters when observations are available at uniformly spaced times by Brockwell et al. \cite{bro4}.
They also generalized the ideas to higher order CARMA processes with nonnegative kernel. The key idea is the decomposition of the CARMA process into a sum of dependent OU processes \cite{bro5}.

Replacing the OU process by a Levy-driven CARMA process with non-negative kernel provides non-negative, heavy-tailed processes with a larger range of auto covariance functions \cite{bro3}.
It is shown that these processes are the convolution of a kernel function with a Levy-driving process.
Gaussian CARMA processes are special cases in which the driving Levy process is Brownian motion. The use of more general Levy processes permits these processes with marginal distributions which may be asymmetric and heavier tailed than Gaussian.
In many situations it is not appropriate to assume Gaussianity of the variables of interest, since the observed time series often exhibit
features like skewness or heavy-tails which contradict the Gaussian assumption.

Jeanblanc et al. \cite{j1} give a representation of self-similar processes with independent increments as stochastic integrals with respect to background driving Levy processes.

Semi-Levy process which is a generalization of Levy process, is an additive process with periodically stationary increments. These processes have been extensively studied by Maejima and Sato \cite{m11}.
Semi-Levy processes are also related to semi-selfsimilar additive processes which has independent increments and are continuous in probability with cadlag paths \cite{m1}.

In this paper we introduce continuous-time ARMA process driven by second order simple semi-Levy measure which has periodically stationary increments. We show that this process is well defined without having to assume further conditions on the driving semi-Levy process. We present the expected value and covariance function of such processes and show that it is associated with a periodically divisible measure. We investigate a new integral representation of such process and discuss on basic properties of these models based on observations made at discrete times.
This study has the potential to provide an approximation for every semi-Levy driven CARMA process.

This paper is organized as follows. In section 2 we present some concepts, theories and ideas regarding the semi-Levy processes, infinitely divisible, self-decomposable distribution and basic properties of them. Section 3 is devoted to the main results and introducing simple semi-Levy driven CARMA processes. For this we present the structure of the measure by a simple semi-Levy compound Poisson measure. We specialize this section to the characteristic function using the concept of periodically divisible measures. We also discuss on the solution of the stochastic differential equation driven by such semi-Levy measure in section 4. For such CARMA process we also study and obtain the first and second moments. The asymptotic behavior and stationarity of the solution are studied in this section. In section 5 we present some simulation of CARMA process driven by simple semi-Levy measure and give an example to illustrate the properties of this process.

\renewcommand{\theequation}{\arabic{section}.\arabic{equation}}
\section{Theoretical framework}
In this section we study the preliminaries such as semi-Levy processes, Levy-Khintchine representation, Levy density and the concepts of infinitely divisible and self-decomposable distributions and their relations to characteristic functions which are used in this paper.

\subsection{Semi-Levy processes}
A general class of stochastic processes with stationary independent increments, called Levy processes are defined in \cite{k1}, \cite{a1}. There are some classical applied probability models which are built on the strength of well-understood path properties of elementary Levy processes. We consider periodic independently scattered random measures, the counterparts of semi-Levy processes in stochastic processes. We provide some basic properties and examples of semi-Levy processes in this section.\\

A stochastic process $\{X_t, t\geq 0\}$ is called an additive process if $X_0= 0$ a.s., it is stochastically continuous, it has independent increments and its sample paths are right-continuous and have left-limits in $t > 0$. Further, if $X_t$ has stationary increments, it is a Levy process. In other words, a more specific definition of Levy process is as following \cite{m1}.

\begin{definition}
A process $\{X_t, t\geq 0\}$ defined on a probability space $(\Omega, \mathcal{F}, \mathbb{P})$ is said to be a Levy process if it possesses the following properties:\\

\item(i) The pathes of $X$ are $\mathbb{P}$-almost surely right continuous with left limits.
\item(ii) $\mathbb{P}(X_0=0)=1$.
\item(iii) For $0\leqslant s\leqslant t$, $X_t-X_s$ is equal in distribution to $X_{t-s}$.
\item(iv) For $0\leqslant s\leqslant t$, $X_t-X_s$ is independent of $\{X_u: u\leqslant s\}$.
\end{definition}
Unless otherwise stated, from now on, when talking of a Levy process, we shall always use the measure $\mathbb{P}$ to be implicity understood as its law.
It is provided a complete characterization of random variables with infinitely divisible distributions via their characteristic functions.
This is the celebrated Levy-Khintchine formula \cite{b01}, \cite{b1}.

\begin{theorem}
The law of Levy process $L$ in terms of triplet $(\gamma, \sigma^2, \nu)$, $\gamma\in\Bbb R$, $\sigma^2\in\Bbb{R^+}$ is infinitely divisible if and only if there exists a triplet such that the characteristic function is given by $E[e^{iuL_t}]=e^{t\psi(u)}$ with
$$\psi(u)=iu\gamma-\sigma^2\frac{u^2}{2}+\int_{-\infty}^{\infty}\big(e^{iux}-1-iux1_{|x|\leq 1}\big)\nu(dx).$$
where the Levy measure $\nu$ satisfies $\nu(\{0\})=0$ and the integrability condition
$\int_{-\infty}^{\infty}1\wedge x^2\nu(dx)<\infty$.
\end{theorem}

As an extension of Levy process, we present the definition and some basic properties and results related of semi-Levy processes \cite{m11}.
\begin{definition}
A subclass of additive processes with the property that for some $p>0$ and for any $s, t\geq 0$,
$$X_{t+p}-X_{s+p}\stackrel{d}{=}X_{t}-X_{s}.$$
where $\stackrel{d}{=}$ denotes the equality in all finite dimensional distributions is called a semi-Levy process with period $p$.
\end{definition}
Linear Brownian motion, compound Poisson and inverse Gaussian processes are familiar processes which are Levy as well.\\

\begin{proposition}
Let $X=\{X_t, t\geq 0\}$ be an additive process and $\mu_t$ be the distribution of $X_t$. If it is a semi-Levy process with period $p$, then
$$\mu_{np+t}=\mu^n_p*\mu_t.$$
for all $n\in{\Bbb N}$ and $t\geq 0$. If the above relation holds for all $n$ and $t\in[0, p)$, then $X$ is a semi-Levy process with period $p$.
\end{proposition}

\subsection{Infinity divisible and self-decomposable distributions}
Let $I({\Bbb R})$ be the class of all infinitely divisible distributions on ${\Bbb R}$.
A distribution $\mu$ is infinitely divisible if for each $n$, there exists a distribution function $\mu_n$ such that $\mu$ is the $n$-fold convolution $\mu_n*\ldots *\mu_n$ of $\mu_n$. The class of possible limit laws consists of the infinitely divisible
distributions.

\begin{remark}
The random variables of a Levy process is infinitely divisible and if $\mu$ an infinitely divisible distribution, we can construct a Levy process from it.
\end{remark}

An important class of random variable models for the unit time distribution as independent effects on the return may need to be scaled to be brought to comparable orders of magnitude before scaling by the square root of $n$ becomes relevant. Such considerations
motivate arbitrary scaling factors and point to self-decomposable laws as candidate models.

\begin{definition}
A probability law of a random variable $X$ is said to self-decomposable just if for every constant $0 <c< 1$ there exists an
independent random variable $X^c$ such that $X\stackrel{d}{=}cX+X^c$.
\end{definition}
The class of self-decomposable distributions, denoted by $L({\Bbb R})$ has the longest history in the study of subclasses of $I({\Bbb R})$. Let $\hat{\mu}(z)$, $z\in {\Bbb R}$ be the characteristic function of $\mu$. Then $\mu$ is said to be self-decomposable if for any
$b>1$, there exists a distribution $\rho_b$ such that
$$\hat{\mu}(z)=\hat{\mu}(b^{-1}z)\hat{\rho}_b(z),$$
where $\rho_b\in I({\Bbb R})$ and $\mu\in L({\Bbb R})$ is also a limiting distribution of normalized partial sums of
independent random variables under infinitesimal condition, and has the stochastic integral representation with respect to a Levy
process.\\

Sato \cite{s1} showed that the Levy (jump) measure of a self-decomposable distribution is always absolutely continuous with respect to the Lebesgue measure and it's density (Levy density) can be characterized by $\nu(x)=\frac{k(x)}{|x|}$, $x\in {{\Bbb R}-{\{0\}}}$
with a so-called $k$-function $k:{\Bbb R}-{\{0\}}\rightarrow {\Bbb R}^+$ which increases on $(-\infty, 0)$ and decreases on $(0,\infty)$.
An infinitely divisible law is self-decomposable if the corresponding Levy density has the above form \cite{c-0}.

\begin{remark}
Self-decomposable laws are infinitely divisible and may be characterized nicely in terms of the Levy density.
\end{remark}
Note that $X(t)$ is a Levy process then $X(1)$ is self-decomposable if and only if $X(t)$ is self-decomposable for every $t>0$.
Levy's continuity theorem enables us to show convergence of distribution through point-wise convergence of characteristic functions.
\begin{theorem}\label{LCT}(Levy's continuity theorem)
If $\varphi_n(u)\rightarrow\varphi(u)$ for every $u$, where $\varphi_n(u)=E[e^{itX_n}]$ and $\varphi$ is continuous at $0$, then $X_n$ converges in distribution to the random variables $X$ with characteristic function $\varphi(u)$.
\end{theorem}

\renewcommand{\theequation}{\arabic{section}.\arabic{equation}}
\section{Semi-Levy driven CARMA process}
If $S=\{S(t), t\in\Bbb R\}$ is a second-order subordinator i.e. nonnegative and nondecreasing Levy process, the semi-Levy driven CARMA$(p,q)$ process $\{Y(t), t\in\Bbb {R^+}\}$, $p>q$ with parameters $a_1, \ldots, a_p, b_0, \ldots, b_q$ is defined via the state space representation of the stochastic differential equation
\be a(D)Y(t)=b(D)DS(t)\label{1}.\ee
where $D$ denotes differentiation with respect to $t$, $a(z)=z^p+a_1z^{p-1}+ \ldots+a_p$, $b(z)=b_0+b_1z+ \ldots+ b_{p-1}z^{p-1}$ and the coefficients $b_j$ satisfy $b_q=1$ and $b_j=0$ for $q<j<p$.
To avoid trivial complications, we shall assume that $a(z)$ and $b(z)$ have no common factors. Since $DS(t)$ does not exist in the usual sense, we interpret the differential equation $(\ref{1})$ by means of its state-space representation, consisting of the observation and state equations
\be Y(t)=\bf b'X(t).\label{2}\ee
and
\be d{\bf X}(t)-{\bf AX}(t)dt={\bf e}dS(t),\label{3}\ee
where $d$ denotes an infinitesimal increment and
\begin{align*}
{\bf A}=\begin{bmatrix}0&1&0&\ldots&0\\
0&0&1&\ldots&0\\
\vdots&\vdots&\vdots&\ddots&\vdots\\
0&0&0&\ldots&1\\
-a_p&-a_{p-1}&-a_{p-2}&\ldots&-a_1\\ \end{bmatrix},
&\hspace{1cm}{\bf e}=\begin{bmatrix}
0\\
0\\
\vdots\\
0\\
1\\
\end{bmatrix},
&\hspace{3mm}{\bf b}=\begin{bmatrix}
b_0\\
b_1\\
\vdots\\
b_{p-2}\\
b_{p-1}\\ \end{bmatrix}.\\
\end{align*}
Every solution of equation $(\ref{3})$ satisfies the following relations for all $t>s$, $s\in\Bbb R$
\be {\bf X}(t)=e^{{\bf A}(t-s)}{\bf X}(s)+\int_s^te^{{\bf A}(t-u)}{\bf e}dS(u).\label{4}\ee
where the integral can be interpreted as the $L^2$-limit of approximating Riemann-Stieltjes sums and also in the path wise sense since the paths of $S$ have bounded variation on compact intervals. From equation $(\ref{4})$ and the independence of the increments of $S$ one can easily verify that
${\bf X}(t)$ is Markov.

\subsection{Structure of semi-Levy measure}
The aim of this section is to present the structure of the simple semi-Levy measure. For this we characterize the measure in Levy-Khintchine representation to an infinitely divisible distribution. This will be done by Levy-Ito decomposition which describes the structure of a general Levy process in terms of three independent auxiliary Levy processes, each with different types of path behavior. In general case, any Levy process may be decomposed into the three independent Levy processes as Brownian motion with drift, compound Poisson process and a square integrable (pure jump) martingale with an a.s. countable number of jumps of magnitude less than 1 on each finite time interval \cite{s1}, \cite{k1}.

\begin{definition}
We call $\{M(0,t], t\geqslant 0\}$, a simple semi-Levy Poisson measure if there exists a partition of the positive real line as $B_i=(s_{i-1}, s_i]$, $i\in\Bbb N$ and $0=s_0< s_1< \ldots$, where for some fixed $r\in\Bbb N$, $|B_i|=|B_{i+kr}|$, $k\in\Bbb N$ and $M(0, \cdot]$ is a Poisson random measure with intensity parameter $\lambda_i$ on $B_i$ for $i \in \Bbb N$, where $\lambda_i=\lambda_{i+kr}$.
\end{definition}

\noindent
Such simple semi-Levy measure has the potential to approximate any semi-Levy measure.
To justify that $M(0,t]$ presented by the above definition is a semi-Levy
measure with period $T$, let $t=kT+s$, $T=\sum_{i=1}^{r}|B_i|$, $s\in B_j$, $j=1, \ldots, r$ and $k=0, 1, \ldots$ then
\vspace{-2mm}
$$M(0, kT+s]=\sum_{i=1}^{kr+j-1}M_i(s_{i-1}, s_i] + M_{kT+j}(s_{j-1}, s].$$
So for all $k\in{\Bbb N}$, $M\big((k-1)T+s, kT+s\big]$ have the same distribution that is the random measure $M$ has periodically stationary increments with period $T$. Thus for $t=(k-1)T+s$, $s\in B_j$, $N(t)=M(0,t]$ is a semi-Levy Poisson process with parameter
\begin{equation}\label{Lam}
\Lambda_t=(k-1)\sum_{i=1}^{r}\lambda_i+\sum_{i=1}^{j-1}\lambda_i+ \frac{\lambda_ja_j^s}{a_j},
\end{equation}
where $r\in\Bbb N$, $a_j=|B_j|$ and $a_j^s=s-s_{j-1}$ for $j=1, \ldots, r$.\\ \\
Let $S(t)$ be a subordinator with the following representation
\begin{equation}\label{LI}
S(t)=\gamma t+ \sum_{k=1}^{N(t)}J_k,
\end{equation}
where $\gamma\in\Bbb R$, $\{N(t), t\geqslant 0\}$ is a simple semi-Levy Poisson process with parameter $\Lambda_t$, defined by (\ref{Lam}) and $\{J_k, k\geqslant 1\}$ is an independent identically distributed sequence of random variables with probability distribution $F$.

\subsection{Representation of periodically divisible random measure}
There is a strong interplay between Levy processes and infinitely divisible distributions and also between semi-Levy processes and periodically divisible (PD) of corresponding random measures. Following the definition of semi-Levy process presented by Sato \cite{m11}, we introduce concept of PD random measure which provides a good platform for obtaining the results of this paper.
\begin{definition}
Random measure $M$ is called periodically divisible with period $T$ if for any fixed $t$, random variables $M(t+(i-1)T, t+iT]$ are independent identically distributed for all indices $i\in{\Bbb Z}$. So the corresponding distribution of $V_k(t)=M(t, t+kT]$ is the same as $k$ times convolution of distribution of $V_1(t)=M(t, t+T]$. Therefore
$$F_{V_k(t)}=F_{V_1(t)}*\ldots*F_{V_1(t)}.$$
\end{definition}
According to Proposition 2.2, we can characterize a PD random measure using its characteristic function.

\begin{definition}
A process $\{X(t), t\geqslant 0\}$ is PD with period $T>0$, if for any $k\in\Bbb N$ and $s\in[0, T)$ we have $X(kT+s)\stackrel{d}{=}X_{T_1}+ X_{T_2}+ \ldots+ X_{T_k}+X^k(s)$, where $X_{T_i}=X(iT)-X((i-1)T)$, $i=1, 2, \ldots, k$ are independent identically distributed random variables which are independent of $X^k(s)=X(kT+s)-X(kT)$. So the corresponding characteristic functions satisfy
$$\varphi_{kT+s}(u)=\big(\varphi_T(u)\big)^k\varphi_s(u),$$
where $\varphi_t(u)=E[e^{iuX(t)}]$.
\end{definition}
It follows that a semi-Levy process is PD with the same period.

\begin{lemma}
Let $\{N(t), t\geqslant 0\}$ be a semi-Levy Poisson process with parameter $\Lambda_t$, defined by (\ref{Lam}) and $\{S_n, n=0, 1, \ldots\}$ be a summation of $n$ independent identically distributed random variables which considered as independent copy of $J$ with characteristic function $\varphi_J(u)=E[e^{iuJ}]$. The jumps $J_j$ are independent of $N(t)$. Then $X(t)=S_{N(t)}$ is called as semi-Levy compound Poisson process with characteristic function.
\begin{equation}\label{Char1}
\varphi_t(u)=e^{\Lambda_t(\varphi_J(u)-1)}.
\end{equation}

\end{lemma}
{\bf Sketch of proof:}
$N(t)$ is a nonhomogeneous process with parameter $\Lambda_t$, so $X$ is stochastically continuous in probability
$$P\{|X(t+h)-X(t)|>\epsilon\}\leq P\{N(t+h)-N(t)>0\}=1-e^{-\Lambda_h}\rightarrow 0,$$
as $h\rightarrow 0$.\\
Thus for any $t\in B_j$, $j\in\Bbb N$ and $t=kT+s$ where $k\in\Bbb N$ and $s\in[0, T)$ we have that
\begin{align}
\varphi_t(u)=E[e^{iuX(t)}]&=\sum_{n=0}^{\infty}P[N(t)=n]E[e^{iuS_n}]=\sum_{n=0}^{\infty}\frac{e^{-\Lambda_t}\Lambda_t^n}{n!}(E[e^{iuJ}])^n\nonumber\\
&=e^{\Lambda_t(E[e^{iuJ}]-1)}=e^{\Lambda_t(\varphi_J(u)-1)}.\nonumber
\end{align}
\begin{example}
The semi-Levy compound Poisson process is an example of PD random measure. Let $X(t)=\sum_{j=1}^{N(t)}J_j$ where $N(t)$ represents the number of jumps and is a semi-Levy Poisson process with parameter $\Lambda_t$. Also $X(t)$ has the characteristic function in the form (\ref{Char1}) with the assumptions of Lemma 3.1. Then we have
\begin{align}
\varphi_t(u)&=E[e^{iuX(kT+s)}]=e^{\Lambda_{kT+s}(\varphi_J(u)-1)}=e^{(k\Lambda_T+\Lambda_s)(\varphi_J(u)-1)}\nonumber\\
&=\big(e^{\Lambda_T(\varphi_J(u)-1)}\big)^ke^{\Lambda_s(\varphi_J(u)-1)}=\big(\varphi_T(u)\big)^k\varphi_s(u),\nonumber
\end{align}
where $\Lambda_T=\sum_{i=1}^{r}\lambda_i$ and $\Lambda_s=\sum_{i=1}^{l-1}\lambda_i+ \frac{\lambda_la_l^s}{a_l}$ where $l<r$.
\end{example}
Since the characteristic function of a random variable determines its distribution, we have a characterization of the distribution of the semi-Levy process by the followings.

\begin{corollary}
By Lemma 3.1, the characteristic function of the semi-Levy process $\{S(t), t\geqslant 0\}$, defined by (\ref{LI}) has the form
\begin{equation}\label{Char2}
E[e^{iuS(t)}]=exp\Big[iu\gamma t+\Lambda_t\big(\varphi_J(u)-1\big)\Big]=e^{iu\gamma t+\Lambda_t \int_{\Bbb R}(e^{iux}-1)\mu(dx)}.
\end{equation}
Since the distribution of $J$ is $\mu$ we have $E[e^{iuJ}-1]=\int_{\Bbb R}(e^{iux}-1)\mu(dx)$ and for some $\gamma\in\Bbb R$, where the semi-Levy measure $\mu$ satisfies $\mu(\{0\})=0$ and $\Lambda_t$ is defined by (\ref{Lam}).
\end{corollary}

\begin{lemma}\label{PD}
If $\{\mu_k\}_{k=0}^{\infty}$ is a sequence of PD random measures with some period $T$ and $\mu_k\rightarrow\mu$, then $\mu$ is also PD random measure with period $T$.
\end{lemma}
{\bf Proof:} The random measure property of $\mu$ follows by a completely similar method to the one of \cite{s1} where such property is proved for infinite divisible distribution and so is valid for corresponding random measures. The periodicity of $\mu$ follows from the fact that all elements of the sequence of measures $\mu_k$ are PD with period $T$ and there is a convergence in distributions to $\mu$ and the period is constant. So the limit of the sequence is also PD with the same period $T$.

\begin{theorem}
Let $S(t)$ be the semi-Levy process defined by (\ref{LI}) and $\gamma\in\Bbb R$, then there exists a simple random measure $\mu$ that it's corresponding characteristic function can be represented by
$$\varphi_t(u)=exp\Big[iu\gamma t+\Lambda_t\int_{-\infty}^{\infty}\big(e^{iux}-1\big)\mu(dx)\Big],$$
where $\Lambda_t$ is defined by (\ref{Lam}).
\end{theorem}
{\bf Proof:} Let $\{\eta_n\}$, $n\in\Bbb N$ be a sequence of real number, monotonic and decreasing to zero and $\{X_n\}$ the sequence has the following characteristic function
$$\varphi_{X_n}(u)=exp\Big[iu\gamma t+\Lambda_t\int_{|x|>\eta_n}\big(e^{iux}-1\big)\mu(dx)\Big].$$
for all $u\in\Bbb R$ and $n\in\Bbb N$, where $\mu$ is the same.
The measure $\mu$ restricted to $\{|x|>\eta_n\}$ is finite and hence $\varphi_{X_n}(u)$ is the characteristic function of a semi-Levy random measure. We clearly have that
$$\lim_{n\rightarrow\infty}\varphi_{X_n}(u)\longrightarrow\varphi_{X}(u),$$
where $\varphi_{X}(u)$ is the characteristic function of $S(t)$ defined by (\ref{Char2}).
By Levy's continuity theorem (Theorem \ref{LCT}) and Lemma \ref{PD}, $\varphi_X$ is the characteristic function of a PD measure provided that $\varphi_X$ is continuous at 0.\\
Continuity of $\varphi_X$ at 0 boils down to the continuity of the integral term. Using the properties of semi-Levy measure and monotone convergence theorem we have
$$|\Lambda_t \int_{\Bbb R}(e^{iux}-1)\mu(dx)|=\Lambda_t \int_{\Bbb R}|e^{iux}-1|\mu(dx)=\Lambda_t \int_{\Bbb R}|\cos ux-1+i\sin ux|\mu(dx)$$
$$\leq \Lambda_t \int_{\Bbb R}\sqrt{(\cos ux-1)^2+\sin^2 ux}\mu(dx)=\sqrt{2}\Lambda_t \int_{\Bbb R}\sqrt{1-\cos ux}\mu(dx)\longrightarrow^{\hspace{-7mm}u\rightarrow 0}0.$$

\section{Characterization of the solution}
Our approach leads to a model which can be interpreted as a solution to the formal differential equation (\ref{1}).
Analyzing the representation of its solution shows that it can be used to define semi-Levy driven CARMA processes.
We take a closer look at the probabilistic properties of the solution, represented in (\ref{4}) such as second moments, Markov
property, non-stationary and limiting distributions and path behavior.
In particular, we characterize the non-stationary distribution and path behavior and give conditions for the existence of it.\\

In order to study the properties of $X$ in (\ref{4}), we need to have the auxiliary results under following conditions.\\ \\
{\bf Condition 1.}
The $X(t)$ is independent of $\{S(r)-S(t), r>t\}$ for all $t\in\Bbb R$.\\ \\
{\bf Condition 2.}
The  eigenvalues of  the matrix {\bf A} are considered to have negative real parts.   We remind that
the eigenvalues of a matrix {\bf A} have negative real parts if and only if
\be \lim_{t\rightarrow\infty}e^{{\bf A}t}=0.\ee
If the above conditions are satisfied the solution (\ref{4}), converges to
\begin{equation}\label{X}
X(t)=\int_{-\infty}^{t}e^{{\bf A}(t-u)}{\bf e}dS(u).
\end{equation}
with the specified properties when $s\rightarrow -\infty$ and $S$ is the semi-Levy process.\\

We have restrict attention to second-order subordinator and for such process represented in (\ref{LI}) and Wald's equation, we have
\begin{equation}\label{ES}
E[S(t)]=\gamma t+\Lambda_t\kappa,
\end{equation}
where $E[J]=\kappa$ and
\begin{equation}\label{VS}
\mbox{var}[S(t)]=\Lambda_t\beta,
\end{equation}
where $E[J^2]=\beta$ and $E[N(t)]=\mbox{var}[N(t)]=\Lambda_t$. Since we have
\begin{align}
\mbox{var}[S(t)]&=E_{N(t)}\big[\mbox{var}\big(S(t)|N(t)\big)\big]+\mbox{var}_{N(t)}\big[E\big(S(t)|N(t)\big)\big]\nonumber\\
&=E_{N(t)}\big[N(t)\mbox{var}(J)\big]+\mbox{var}_{N(t)}\big[\gamma t+N(t)E[J]\big]\nonumber\\
&=\Lambda_t\big(\mbox{var}(J)+E^2[J]\big)=\Lambda_t\beta.\nonumber
\end{align}
In the following we find mean and second moment of $X(t)$, represented in (\ref{X}), for $t\in B_{kr+j}$. This is by the fact as we have a partition on positive real line, $B_i=(s_{i-1}, s_i]$, $i\in\Bbb N$, described more in Definition 3.1.\\ \\
Let $\Pi=\{0=s_0< s_1< \ldots< s_{r-1}< s_{r}=T< s_{r+1}< \ldots< s_{kr}< \ldots< s_{kr+j-1}\leqslant s_{kr+j}=t\}$
and $t=kT+s$ be a point inside some interval, say $B_{kr+j}$. For the simplicity we can re-consider this point as the last point of the last period interval. By this we have $r$ subinterval in each period interval. So this new notation implies that
\begin{align}
X(t)&=\int_{-\infty}^{kT+s}e^{{\bf A}(t-u)}{\bf e}dS(u)
=\sum_{n=1}^{\infty}\int_{(k-n)T+s}^{(k-n+1)T+s}e^{{\bf A}(t-u)}{\bf e}dS(u)\nonumber\\
&=\lim_{d\rightarrow\infty}\sum_{n=0}^{d}\sum_{i=0}^{r-1}e^{{\bf A}t}\int_{(k-n)T+s_{j-2-i}}^{(k-n)T+s_{j-1-i}}e^{-{\bf A}u}{\bf e}dS(u)
+e^{{\bf A}t}{\bf e}\int_{kT+s_{j-1}}^{kT+s}e^{-{\bf A}u}dS(u).\nonumber
\end{align}
So the expected value of $X(t)$ for large $r$ is
\begin{align}
E[X(t)]=E\Big[\int_{-\infty}^te^{{\bf A}(t-s)}{\bf e}dS(s)\Big]
&=\lim_{d\rightarrow\infty}\sum_{n=0}^{d}\sum_{i=0}^{r-1}e^{{\bf A}t}\int_{(k-n)T+s_{j-2-i}}^{(k-n)T+s_{j-1-i}}e^{-{\bf A}u}{\bf e}E[dS(u)]\nonumber\\
&+e^{{\bf A}t}\int_{kT+s_{j-1}}^{kT+s}e^{-{\bf A}u}{\bf e}E[dS(u)].\nonumber
\end{align}
For $j$ outside the range $1, \ldots, r$ we assume that $s_j$ is equal to $s_{j+r}$.
By the fact that $S(t)$ is semi-Levy process with periodically stationary increments and using (4.11) and (3.5) we have
$$E[S(s_{i+1})]-E[S(s_{i})]=\gamma\Delta s_i+\kappa(\Lambda_{s_{i+1}}-\Lambda_{s_{i}}),$$
where $\Delta s_i=s_{i+1}-s_{i}$ and we have $\Lambda_{s_{i+1}}-\Lambda_{s_{i}}=\frac{\lambda_i}{a_i}$. For any $u, u+du\in B_i$, $i=0, \ldots, r-1$
$$E[dS(u)]=E[S(u+du)]-E[S(u)]=(\gamma+\frac{\lambda_i}{a_i}\kappa)du.$$
Therefore
\begin{align}
&E[X(t)]=\lim_{d\rightarrow\infty}\sum_{n=0}^{d}\sum_{i=0}^{r-1}e^{{\bf A}t}(\gamma+\frac{\lambda_i}{a_i}\kappa)
\int_{(k-n)T+s_{j-2-i}}^{(k-n)T+s_{j-1-i}}e^{-{\bf A}u}{\bf e}du
+e^{{\bf A}t}(\gamma+\frac{\lambda_i}{a_i}\kappa)\int_{kT+s_{j-1}}^{kT+s}e^{-{\bf A}u}{\bf e}du\nonumber\\
&=\lim_{d\rightarrow\infty}\sum_{n=0}^{d}\sum_{i=0}^{r-1}e^{{\bf A}(kT+s)}(\gamma+\frac{\lambda_i}{a_i}\kappa)
(-{\bf A}^{-1})\big(e^{-{\bf A}[(k-n)T+s_{j-1-i}]}-e^{-{\bf A}[(k-n)T+s_{j-2-i}]}\big){\bf e}\nonumber\\
&+e^{{\bf A}(kT+s)}(\gamma+\frac{\lambda_j}{a_j}\kappa)(-{\bf A}^{-1})\big(e^{-{\bf A}[kT+s]}-e^{-{\bf A}[kT+s_{j-1}]}\big){\bf e}.\nonumber\\
&=-{\bf A}^{-1}e^{{\bf A}s}\sum_{n=0}^{\infty}e^{{\bf A}nT}\sum_{i=0}^{r-1}(\gamma+\frac{\lambda_i}{a_i}\kappa)\big(e^{-{\bf A}s_{j-1-i}}-e^{-{\bf A}s_{j-2-i}}\big){\bf e}\nonumber\\
&-{\bf A}^{-1}e^{{\bf A}s}(\gamma+\frac{\lambda_j}{a_j}\kappa)\big(e^{-{\bf A}s}-e^{-{\bf A}s_{j-1}}\big){\bf e}\nonumber\\
&=-{\bf A}^{-1}e^{{\bf A}s}\Big[\frac{I}{I-e^{{\bf A}T}}\sum_{i=0}^{r-1}(\gamma+\frac{\lambda_i}{a_i}\kappa)\big(e^{-{\bf A}s_{j-1-i}}-e^{-{\bf A}s_{j-2-i}}\big)+(\gamma+\frac{\lambda_j}{a_j}\kappa)\big(e^{-{\bf A}s}-e^{-{\bf A}s_{j-1}}\big)\Big]{\bf e}.\nonumber
\end{align}
where $I$ is the identity matrix. So we arrive at the following lemma.

\begin{lemma}
The expected value of $X(t)$ defined by (\ref{X}) is periodic with period $T$, that is $E[X(t+T)]=E[X(t)]$.
So for $t=kT+s$, the solution $Y(t)$ of the ordinary differential equation is also periodic.
\end{lemma}
Now we find the covariance function of $X(t)$ as $\mbox{cov}\big(X(t), X(t+h)\Big)$ for $h>0$ and $t=kT+s$ and $t+h=k'T+s'$. Let $\tilde{X}(t)=X(t)-E[X(t)]$, so
$$\mbox{cov}\big(X(t), X(t+h)\Big)=E[\tilde{X}(t)\tilde{X}'(t+h)].$$
Using the partition $\Pi$ and by the fact that the increments of semi-Levy process $S(t)$ are periodically stationary and independent, for $s_i<s'_i$ and $u<z$ we have
\begin{align}
E[\tilde{X}(t)\tilde{X}'(t+h)]&=\lim_{d\rightarrow\infty}\sum_{n=0}^{d}\sum_{i=0}^{r-1}e^{{\bf A}t}
\Big[\int_{(k-n)T+s_{j-i-2}}^{(k-n)T+s_{j-i-1}}e^{-{\bf A}u}{\bf e}{\bf e'}e^{-{\bf A'}u}\mbox{var}[dS(u)]\Big]e^{{\bf A'}(t+h)}\nonumber\\
&+e^{{\bf A}t}\Big[\int_{kT+s_{j-1}}^{kT+s}e^{-{\bf A}u}{\bf e}{\bf e'}e^{-{\bf A'}u}\mbox{var}[dS(u)]\Big]e^{{\bf A'}(t+h)}.\nonumber
\end{align}
and using the variance of $S(t)$ as $\mbox{var}[S(t)]=\beta\Lambda_t$ and for any $u, u+du\in B_i$, $i=0, \ldots, r-1$
\begin{align}
\mbox{var}[S(u+du)-S(u)]&=\mbox{var}[S(u+du)]+\mbox{var}[S(u)]-2\mbox{cov}\big(S(u+du), S(u)\big)\nonumber\\
&=\mbox{var}[S(u+du)]-\mbox{var}[S(u)]=\beta(\Lambda_{u+du}-\Lambda_{u})=\beta\frac{\lambda_i}{a_i}du.\nonumber
\end{align}
This is by the fact that independent increments of $S$ implies that $\mbox{cov}\big(S(u+du), S(u)\big)=\mbox{var}[S(u)]$. Therefore
\begin{align}
&E[\tilde{X}(t)\tilde{X}'(t+h)]=\beta e^{{\bf A}(kT+s)}\Big[\sum_{n=0}^{\infty}\sum_{i=0}^{r-1}\frac{\lambda_i}{a_i}
\int_{(k-n)T+s_{j-i-2}}^{(k-n)T+s_{j-i-1}}e^{-{\bf A}u}{\bf e}{\bf e'}e^{-{\bf A'}u}du\nonumber\\
&+\frac{\lambda_j}{a_j}\int_{kT+s_{j-1}}^{kT+s}e^{-{\bf A}u}{\bf e}{\bf e'}e^{-{\bf A'}u}du\Big]e^{{\bf A'}(kT+s)}e^{{\bf A'}h}\nonumber\\
&=\beta e^{{\bf A}kT}e^{{\bf A}s}\Big[\sum_{n=0}^{\infty}\sum_{i=0}^{r-1}\frac{\lambda_i}{a_i}
\int_{s_{j-i-2}}^{s_{j-i-1}}e^{-{\bf A}[(k-n)T+u]}{\bf e}{\bf e'}e^{-{\bf A'}[(k-n)T+u]}du\nonumber\\
&+\frac{\lambda_j}{a_j}\int_{s_{j-1}}^{s}e^{-{\bf A}[kT+u]}{\bf e}{\bf e'}e^{-{\bf A'}[kT+u]}du\Big]e^{{\bf A'}kT}e^{{\bf A'}s}e^{{\bf A'}h}\nonumber\\
&=\beta e^{{\bf A}s}\Big[\sum_{n=0}^{\infty}e^{{\bf A}nT}\Big(\sum_{i=0}^{r-1}\frac{\lambda_i}{a_i}
\int_{s_{j-i-2}}^{s_{j-i-1}}e^{-{\bf A}u}{\bf e}{\bf e'}e^{-{\bf A'}u}du\Big)e^{{\bf A'}nT}
+\frac{\lambda_j}{a_j}\int_{s_{j-1}}^{s}e^{-{\bf A}u}{\bf e}{\bf e'}e^{-{\bf A'}u}du\Big]e^{{\bf A'}s}e^{{\bf A'}h}.\nonumber
\end{align}
So we arrive at the following lemma.\\
\begin{lemma}
The covariance function of $X(t)$ defined by (\ref{X}) is periodic with period $T$. In the other hand
$$\mbox{cov}\big(X(t+T), X(t+h+T)\Big)=\mbox{cov}\big(X(t), X(t+h)\Big).$$
So for $t=kT+s$, $t+h=k'T+s'$ and Lemma 4.1, the solution $Y(t)$ of the ordinary differential equation is weakly periodically stationary.\\
\end{lemma}

\begin{remark}
The eigenvalues of a matrix {\bf A}, denoted by $\lambda_1, \ldots, \lambda_p$ are the same as the zeroes of the autoregressive polynomial $a(z)$.
\end{remark}

\begin{proposition}
If $S$ is a second-order subordinator and the eigenvalues of a matrix {\bf A} have negative real parts, then the semi-Levy driven CARMA$(p,q)$ process with parameters $a_1, \ldots, a_p, b_0, \ldots, b_q$ by equations (\ref{2}) and (\ref{X}) is defined as
\be Y(t)={\bf b'X(t)}=\int_{-\infty}^{\infty}h(t-u)dS(u),\label{Y}\ee
where the function $h(t)={\bf b'}e^{{\bf A}t}{\bf e}I_{(-\infty, t]}$ is called the kernel of the CARMA process $\{Y(t)\}$. When the Condition 1 is satisfied, then $Y$ is a causal function of $S$.
\end{proposition}

\begin{remark}
The representation (\ref{Y}) shows that if the kernel $h$ is nonnegative then, the semi-Levy driven CARMA process is nonnegative and can be used to represent nonnegative quantity process such as stochastic volatility.
\end{remark}

\section{Simulation}
In this section, we verify the theoretical results concerning periodically stationary  structure of the output of the CARMA model (3.2) and (3.3)  by simulation.
For this, we consider a discretization of the process by imposing a partition for the whole duration of it with equally spaced points of length one. Then we generate the discretized  version of the CARMA process $Y(t)$ by the following procedure which provides a discrete time periodically stationary  process.
 In the followings we briefly describes the simulation steps.\\

We  present a simulation method for some discretized version of   CARMA processes driven by some simple semi-Levy measure.
Following Definition 3.1  we  simulate some simple semi-Levy Poisson measure.
Let $T>0$ be the period of the increments of the measure and consider  period intervals as $\big((n-1)T, nT\big], \;\; n \in \Bbb N$.
We assume that the number of partitions in each period is fixed, say $r$,  and the subintervals $B_1, B_2, \ldots , B_r$ form a partition of first period interval.
The following period intervals are partitioned by the same number and same lengths of subintervals.
Successive elements of the partitions of successive scale intervals are denoted by $B_1, B_2, \ldots $  where their lengths admit
equalities $|B_{i+kr}|=|B_i|, \;\; i,k \in \Bbb N$.  The simple semi-Levy process $S(t)$, defined by (3.6), on  partition $B_i$  has Poisson distribution with intensity parameter $\lambda_i$ where $\lambda_{i+kr}=\lambda_i$ for  $i \in \Bbb N$.
Therefore, the simulation algorithm of $S(t)$ is determined as follows.

\begin{itemize}
\item[1.] Consider some positive value $T$ as the length of the period of {simple semi-Levy process} $S(t)$, defined by (3.6), and some integer $M$ as the number of corresponding period intervals for simulation.

\item[2.] Decide about the number of elements of partition in each period, say $r$.

\item[3.] Consider $r$ different positive real numbers $l_1, l_2, \cdots, l_r$ so that $T=\sum_{i=1}^r l_i$ and a  partition of first period interval $(0, T]$ by $B_1, B_2, \cdots, B_r$ where $B_i=(s_{i-1}, s_i]$, $s_0=0$ and $s_i=\sum_{j=1}^i l_j$ for  $i=1,2, \cdots, r$.
    Elements of partitions of successive period intervals $\big((i-1)T, iT]$ for $i=2, \cdots, M$ are $B_{r+1}, B_{r+2}, \cdots, B_{rM}$ where
    $\; |B_{i+kr}|=|B_i|=l_i, \; i=1, 2, \cdots, rM$.

\item[4.] Let the positive real numbers $\lambda_1, \lambda_2, \cdots, \lambda_r$ be chosen as Poisson rates of occurrences, corresponding to the increments of $N(t)$,  in (3.6), on  $B_1,B_2, \cdots , B_r$.
    Also increments of $N(t)$ on $B_{i+kr}$ has Poisson rate of occurrence $\; \lambda_{i+kr}=\lambda_i, \; i=r+1, r+2,\cdots , rM$.

\item[5.] Generate an independent sequence of Poisson random variables $U_i$ with parameter $\lambda_il_i$ on  $B_i$  for $\; i=1,2,\cdots , rM$ as
 $u_1, u_2, \cdots , u_{rM}$.

\item[6.] Create independent samples $u_i$ from Uniform distribution $U(s_{i-1}, s_i]$.
Then sort these samples on each $B_i$  and denote these ordered samples by  $t_{i,1}, \cdots , t_{i,u_i}$  for $\; i=1,2,\cdots , rM$. Finally  evaluate  $N(t)= \sum_{i=1}^{rM} \sum_{j:t_{i,j}\leq t} 1$ for $t \in (0, MT]$\footnote{ This is by the fact that sample points of occurrence  of Poisson random variables on an interval  follows the order statistics of Uniform distribution}.

\item[7.] Let $\gamma $ be a real number and generate independent identically distributed random variables $J_i,  \; i=1, 2, \cdots , rM$ from some probability distribution $F$, say standard Normal distribution. Then determine $S(t)$ by relation (3.6).
\end{itemize}

After evaluating $S(t)$, we produce the CARMA process in (3.2) and (3.3)  by the following steps.
\begin{itemize}
\item[1.] Conside $p$ as some specified integer value.

\item[2.]  Following   condition 2 and Remark 4.1
consider $p$  roots for the autoregressive polynomial $a(z)$  with negative real parts and calculate
  the coefficients of $a(z)$ as  $a_1, a_2, \ldots, a_p$.

\item[3.] Create the matrix $A$.

\item[4.] Consider some real  values for  the parameters $b_0, b_1, \ldots, b_{q-1}$, so that $a(z)$ and $b(z)$ have no common factors.

\item[5.] Provide a discretization of $X(t)$ in (3.3) by imposing an equally spaced  partition for $(0, MT]$ with some small
space $h>0$. Also consider $p$ initial values of such discretized vector ${\bf X}(\cdot )$.

\item[6.] Finally, using the discretized value of $X(t)$ provided by previous step and the values considered for  parameters $b_0, b_1, \ldots, b_{q-1}$ and relation (3.2), evaluate the corresponding discretized values of $Y(t)$.
\end{itemize}

We simulate the process  using the proposed algorithm. Hurd and Gerr \cite{h-0} presented the graphical methods to verify that the simulated series are indeed PC.
Soltani and Azimmohseni \cite{s2} and Hurd and Miamee \cite{h-1}(Chapter 10) used the same diagnostic method to check whether the simulated data are PC with period $T$.\\
In this method, for a sample of size $n$, $X(0), \ldots, X(n-1)$ and a fixed $\bf M$, they
plotted the significant values of the sample spectral coherence
$$
|\gamma(p,q,{\bf M})|^2=\frac{|\sum_{m=0}^{{\bf M}-1}d_X(\theta_{p+m})\bar d_X(\theta_{q+m})|^2}{\sum_{m=0}^{{\bf M}-1}|d_X(\theta_{p+m})|^2 \sum_{m=0}^{{\bf M}-1}|d_X(\theta_{q+m})|^2}.
$$
against $p, q=0,\ldots, n-1$. It has the non-zero values for $|p - q|=cn/T$, $c=0,\ldots, T-1$, where
$$d_X(\theta_{p})=\sum_{n=0}^{n-1}X(t)e^{it\theta_{p}},\hspace{.5cm}\theta_{p}=\frac{2\pi p}{n}, \hspace{.5cm}p=0,1,,\ldots,n-1. $$
Therefore, we can say something about the nature of the analyzed time series: 
\begin{itemize}
\item If in the square only the main diagonal appears, so $X(t)$ is a stationary time series.
\item If there are some significant values of statistic and they seem to lie along the parallel equally spaced diagonal lines, then $X_t$ is likely PC-T, where $T$ is the “fundamental” line spacing. Algebraically, $T$ would be the gcd of the line spacings from the diagonal; for a sequence to be PC-T, not all lines are required to be present.
\item If there are some significant values of statistic but they occur in some non-regular places, then $X(t)$ is a nonstationary time series in other than periodic sense; but note there are many hypotheses being tested, so some threshold exceedances are to be expected.
\end{itemize} 
In the following, we provide one example to investigate our process.\\\\
\noindent
{\bf Example:} We provide a simulated discretized semi-Levy driven CARMA process. For this we consider the parameters of the simple semi-Levy process $S(t)$ as $r=7$, $k=40$ and the length of successive subintervals of each period intervals as $2,2,2,2,1,1,2$ where corresponding rate of Poisson occurrence on these subintervals are assumed as $6,4,2,10,4,8,12$. Also the random variables  $J_k$ is considered to have  Normal distribution with mean 3 and variance 1.\\
In this example we  simulate a  CARMA(3,2) model by assuming the roots of  $a(z)$ as $z_1=-1$, $z_2=-2-i$ and $z_3=-2+i$.  So the value of parameters $a_1$, $a_2$ and $a_3$ are determined as  $5,9$ and $5$ respectively. Thus, the matrix A is
\begin{align*}
{\bf A}=\begin{bmatrix}
0&1&0\\
0&0&1\\
-5&-9&-5\\
\end{bmatrix},
\end{align*}
By the CARMA(3.2) model in the form (3.3),  we generate  the values of discretized  $X(t)$ by considering some equally space partition for the duration of $M=40$ period intervals. Finally, we consider the values of parameters as $b_0=0.5$, $b_1=2$ and $b_2=1$ and produce the CARMA process. 	
Then we follow to verify the output of the model which provides a periodically stationary process.\\
In Figure 1, we see the simulated data of size $n=480$ with the suggested simulation algorithm (top).
The sample autocorrelation plot of this process (bottom left) and
 the sample coherent statistic $|\gamma(p,q,{\bf M})|^2$ of data (bottom right).
The parallel lines for the sample spectral coherence confirm that
the simulated data are PC.  Also, In this plot, the first significant off-diagonal is at $|p-q| = 40 $ which verifies the first significant peak at 40 and hhis shows that there is a second-order PC structure with period $T = 480/40 = 12$ in the data.

\begin{figure}[h]
\begin{center}
 \includegraphics[width=16cm, height=9.5cm,  angle=0] {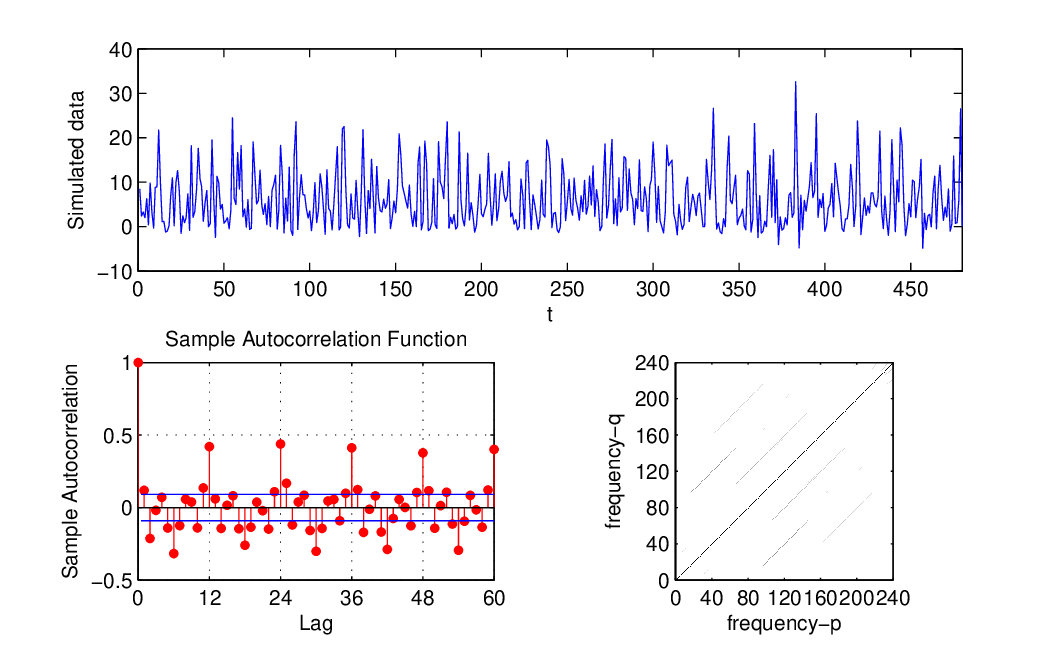}
\end{center}
\vspace{-.81cm}
\caption{\footnotesize{Top: the simulated data of size $n=480$; bottom left: the sample autocorrelation plot of this process; bottom right: the significant values of the
sample spectral coherence with $\alpha=0.01$ and $\mathbf M=40$}.\label{Fig1}}
\end{figure}

\end{document}